\documentclass[12pt,twoside]{article}
\usepackage{amssymb,indentfirst}

{\catcode`\@=11
\gdef\ps@myheadings{\let\@mkboth\@gobbletwo
 \def\@oddhead
          {\vbox{\noindent\hfill
                {\small\it Local subexponential behaviour}
                                         \hfill\rm\thepage\vskip 4pt\hrule}}%
 \def\@oddfoot{}
 \def\@evenhead
          {\vbox{\noindent\rm\thepage\hfil\small\it
          S.\ Asmussen, S. Foss and D. Korshunov\hfil\vskip 4pt \hrule}}%
 \def\@evenfoot{}\def\chaptermark##1{}\def\sectionmark##1{}%
 \def\subsectionmark##1{}}

\gdef\@begintheorem#1#2{\trivlist \item[\hskip \labelsep{\indent\bf #1\ #2.}]\it}
\gdef\@opargbegintheorem#1#2#3{\trivlist
      \item[\hskip \labelsep{\bf #1\ #2\ (#3)}]\it}
\gdef\@endtheorem{\endtrivlist}
 }

\font\scal=rsfs10 at 12pt

\newcommand{\calS}{\mbox{\scal S}}
\newcommand{\calL}{\mbox{\scal L}}
\newcounter{remark}
\newcommand\remark
       {\refstepcounter{remark}
                  {R\,e\,m\,a\,r\,k~\arabic{remark}.\ }}
\newcommand\proof{P\,r\,o\,o\,f}
\newcounter{definition}
\newcommand\definition
       {\refstepcounter{definition}
                  {\bf Definition~\arabic{definition}.\ }}

\newlength{\myskip}
\newtheorem{Theorem}{Theorem}
\newtheorem{Proposition}{Proposition}
\newtheorem{Lemma}{Lemma}
\newtheorem{Corollary}{Corollary}

\newcommand\mysection[1]{
              \refstepcounter{section}
              \section*{\centering\normalsize\bf\thesection.~#1}
                         }

\renewcommand{\P}{{\bf P}}
\newcommand{\E}{{\bf E}}

\begin{document}
\thispagestyle{empty}

\section*{\centering\large\bf
           Asymptotics for sums
of random variables\\
with local subexponential behaviour\footnote{Partially
supported by INTAS grant No.~265,
RFBR grant No.~99-01-01561,
and EPSRC grant No.~R58765/01}}

\section*{\centering\normalsize\slshape\bfseries
S\o ren Asmussen\footnote{Department of Mathematical Statistics,
Lund  University, Box 118, 221 00 Lund, Sweden. E-mail:
asmus@maths.lth.se},
Serguei Foss and Dmitry
Korshunov\footnote{Sobolev Institute of Mathematics,
630090 Novosibirsk, Russia
and Heriot-Watt University, Edinburgh, EH14 4AS, UK.
E-mail: S.Foss@ma.hw.ac.uk, D.Korshunov@ma.hw.ac.uk}}

\begin{abstract}\noindent
We study distributions $F$ on $[0,\infty)$ such that
for some $T\le\infty$, $F^{*2}(x,x+T]\sim 2 F(x,x+T]$.
The case $T=\infty$ corresponds to $F$ being subexponential,
and our analysis shows that the properties for $T<\infty$ are,
in fact, very similar to this classical case.
A parallel theory is developed in the presence of densities.
Applications are given to random walks, the key renewal theorem,
compound Poisson process and Bellman-Harris branching processes.
\\[2mm]
{\it Keywords}: sums of independent random variables,
subexponential distributions, distribution tails,
local probabilities.
\end{abstract}

\mysection{Introduction}\label{Introduction}

For a probability distribution $F$ on the real line,
let $F(x)=F(-\infty,x]$ denote the
distribution function and $\overline F(x)=F(x,\infty) = 1 -F(x)$
the tail. The class $\calS$ of subexponential distributions
is defined by the requirement
$\overline{F^{*2}}(x)\sim 2\overline F(x)$ as $x\to\infty$
($F^{*n}=$ $n$th convolution power) and that
the support is contained in $[0,\infty)$. This class plays an
important role in many applications
(see, e.g., [\ref{Ch}, \ref{EG}, \ref{Ver},
\ref{RSST}, \ref{A} Ch.\ IX].
For example, one of the key results in the theory is:

\begin{Theorem}\label{fundament}
Let  $S_n=\xi_1+\cdots+\xi_n$ be a sequence of
partial sums of i.i.d.\ random variables with  common
distribution $F$, and let $\tau$ be an independent
integer-valued random variable. If $F\in\calS$ and
${\bf E}(1+\delta)^\tau<\infty$ for some $\delta>0$, then
${\bf P}(S_\tau >x) \sim {\bf E}\tau\cdot\overline F(x)$
as $x\to\infty$.
\end{Theorem}
Special cases of this result provide asymptotics for tails
of waiting times in the GI/G/1 queue, for ruin probabilities
and Bellmann--Harris branching processes (see further the
references in later parts of this paper).

Any subexponential distribution is long-tailed, i.e., for any fixed $T$,
$\overline F(x+T)\sim \overline F(x)$ as $x\to\infty$.
This easily yields $F^{*n}(x,x+T]$ $=o(\overline F(x))$ for all
$T<\infty$ and all $n$.
Some applications, however, call for more detailed properties
of $F^{*n}(x,x+T]$ when $T<\infty$, but the theory is more
scattered so the references that we know of are few:
Chover {\em et al.}
[\ref{CNW}, Section 2] gave local theorems for some classes of lattice
distributions; densities were considered in [\ref{CNW}, Section 2]
(requiring continuity) and in Kl\"uppelberg
[\ref{Claudia}] who considered asymptotics of densities
for a special case (see also Sgibnev [\ref{S}]
for some results on the densities on ${\bf R}$); and finally
Bertoin and Doney [\ref{BD}] and
Asmussen {\em et al.} [\ref{AKKKT}]
dealt with the case where $F$ is the ladder height
distribution in a random walk in order to provide more
detailed asymptotics of the random walk maximum than
the standard consequences of Theorem \ref{fundament}.

The aim of the present paper is to develop a more systematic theory.
Fix $0<T\le\infty $ and write $\Delta=(0,T]$,
\begin{eqnarray*}
x+\Delta &\equiv& \{x+y:y\in\Delta\}=(x,x+T],
\quad x\in{\bf R}.
\end{eqnarray*}
Motivated from  [\ref{AKKKT}], we call $F$
(concentrated on $[0,\infty)$) {\em $\Delta$-subexponential}
if the function $F(x+\Delta)$ is long-tailed
(see Definition \ref{def.long.tailed} below)
and $F^{*2}(x+\Delta)\sim 2F(x+\Delta)$ (where $g(x)\sim h(x)$
means that $g(x)/h(x)\to 1$, $x\to\infty$).
Here $T=\infty$ corresponds to ordinary subexponential distributions.
We will see that all standard examples of subexponential
distributions are also $\Delta$-subexponential when $T<\infty$,
and that the standard theory for $T=\infty$
carries over to $T<\infty$ practically without changes.
We thereby provide a general theory covering both the classical
subexponential case and some of the more refined questions
encountered in [\ref{AKKKT}], and we also give some further
applications motivating this generalization, see for example
the results from renewal theory in Section \ref{KRT}.

In Section \ref{Delta}, we  derive the properties
of $\Delta$-subexponential distributions and prove
a natural analogue of Theorem \ref{fundament}.
In Section \ref{Local}, we define distributions with
subexponential densities and study their properties.
In Section \ref{Sufficient}, sufficient conditions
for $\Delta$-subexponentiality are given.
In Section \ref{Supremum}, we apply results from Sections
\ref{Delta} and \ref{Local} to the asymptotic description
of the distribution of the supremum of a random walk with
negative drift. The rest of the paper contains
further applications to Compound
Poisson Processes, Infinitely Divisible Laws,
Bellman--Harris Branching Processes and the Key Renewal Theorem.

\mysection{$\Delta$-Subexponential distributions}
\label{Delta}

\definition\label{def.long.tailed}
We say that a distribution $F$ on ${\bf R}$
belongs to the class $\calL_\Delta$
if $F(x+\Delta)>0$ for all sufficiently
large $x$ and
\begin{equation}\label{lt}
\frac{F(x+t+\Delta)}{F(x+\Delta)} \to 1
\quad \mbox{as} \quad x\to\infty,
\end{equation}
uniformly in $t\in [0,1]$.

Calling a function $g(x)$ {\em long-tailed} if
$g(x+t)/g(x)\to 1$ uniformly in $t\in[0,1]$, we see that
the definition is equivalent to $F(x+\Delta)$ being long-tailed.
If $T=\infty$, then we write $\calL$ instead of $\calL_\Delta$  and say that
$F$ is long-tailed.
It follows from the definition that  one can choose
a function $h(x)\to\infty$ such that (\ref{lt})
holds uniformly in $|t|\le h(x)$.

\begin{Proposition}\label{long.tailed}
Let the distributions $F$ and $G$ belong to the class
$\calL_\Delta$ for some $\Delta$.
Then $F*G\in\calL_\Delta$ and
\begin{eqnarray}\label{long.tailed.below}
\liminf_{x\to\infty}
\frac{(F*G)(x+\Delta)}{F(x+\Delta)+G(x+\Delta)}
&\ge& 1.
\end{eqnarray}
\end{Proposition}

\proof. Let $\xi$ and $\eta$ be two independent random
variables with corresponding distributions $F$ and $G$.
Take an increasing function $h(x)\uparrow\infty$ such
that $h(x)<x/2$, $F(x-y+\Delta)\sim F(x+\Delta)$ and
$G(x-y+\Delta)\sim G(x+\Delta)$ as $x\to\infty$ uniformly
in $|y|\le h(x)$. Consider the event
$B(x,t)=\{\xi+\eta\in x+t+\Delta\}$. The estimate
(\ref{long.tailed.below}) follows from the inequality
\begin{eqnarray*}
\P(B(x,0)) &\ge&
\P(B(x,0),\ |\xi|\le h(x))+\P(B(x,0),\ |\eta|\le h(x))
\end{eqnarray*}
combined with
\begin{eqnarray*}
\P (B(x,0),\ |\xi|\le h(x))
&=& \int_{-h(x)-0}^{h(x)} G(x{-}y{+}\Delta)F(dy)\\
&\sim& G(x{+}\Delta)\int_{-h(x)-0}^{h(x)} F(dy)
\sim G(x{+}\Delta),\\
\P (B(x,0),\ |\eta|\le h(x)) &\sim& F(x{+}\Delta).
\end{eqnarray*}

The probability of the event $B(x,t)$ is equal to the sum
\begin{eqnarray*}
\lefteqn{\P(B(x,t),\ \xi\le x-h(x))+\P(B(x,t),\ \eta\le h(x))}\\
&& \hspace{15mm} +\P(B(x,t),\ \xi>x-h(x),\ \eta>h(x))
\ \equiv\ P_1(x,t)+P_2(x,t)+P_3(x,t).
\end{eqnarray*}
In order to prove that $F*G\in\calL_\Delta$,
we need to check that ${\bf P}(B(x,t))\sim{\bf P}(B(x,0))$
as $x\to\infty$ uniformly in $t\in[0,1]$. This
follows from the relations
\begin{eqnarray*}
P_1(x,t) &=&
\int_{-\infty}^{x-h(x)} G(x+t-y+\Delta)F(dy)\\
&\sim& \int_{-\infty}^{x-h(x)} G(x-y+\Delta)F(dy)= P_1(x,0),
\end{eqnarray*}
$P_2(x,t) \sim P_2(x,0)$, by the same reasons, and
\begin{eqnarray*}
P_3(x,t) &=& \int_{x-h(x)}^{x-h(x)+t+T}
{\bf P}(\eta\in x+t-y+\Delta,\eta>h(x)) F(dy)\\
&\le& {\bf P}(\eta>h(x)) F(x-h(x)+(0,t+T])
=o(F(x+\Delta)).
\end{eqnarray*}

By induction, Proposition
\ref{long.tailed} yields

\begin{Corollary}\label{long.tailed.n}
Let $F\in\calL_\Delta$ for some $\Delta$.
Then, for any $n\ge2$, $F^{*n}\in\calL_\Delta$ and
\begin{eqnarray*}
\liminf_{x\to\infty}\frac{F^{*n}(x+\Delta)}{F(x+\Delta)}
&\ge& n.
\end{eqnarray*}
\end{Corollary}

\definition Let $F$ be a distribution on ${\bf R}^+$
with unbounded support. We say that $F$ is
{\it $\Delta$-subexponential} and write $F\in\calS_\Delta$ if
$F\in\calL_\Delta$ and
\begin{eqnarray*}
(F*F)(x+\Delta) &\sim&
2 F(x+\Delta)\quad\mbox{as }x\to\infty.
\end{eqnarray*}

Equivalently, a random variable $\xi$ has a
$\Delta$-subexponential distribution if the function
$\P(\xi\in x+\Delta)$ is long-tailed and,
for two independent copies $\xi_1$ and $\xi_2$ of $\xi$,
\begin{eqnarray*}
\P (\xi_1+\xi_2\in x+\Delta) &\sim&
2\P(\xi\in x+\Delta) \quad \mbox{as }x\to\infty.
\end{eqnarray*}

\remark The class of ${\bf R}^+$-subexponential
distributions coincides with the standard class
$\calS$ of subexponential distributions.
Typical examples of $\calS_\Delta$
distributions (for all $T>0$) are the same,
in particular
the Pareto--, lognormal--, and Weibull
(with parameter between 0 and 1) distributions,
as will be shown in Section \ref{Sufficient}.
Also, many properties of
$\calS_\Delta$-distributions with finite $\Delta$
are very close to those of subexponential
distributions, as will be shown below.
However, a main difference is that for $T<\infty$,
the function $F(x+\Delta )$ may be non-monotone in $x$,
whereas it is non-increasing for $T=\infty$.

\remark It follows from the definition that, if $F\in\calS_\Delta$
for some finite interval $\Delta = (0,T]$, then $F\in\calS_{n\Delta}$
for any $n=2,3,\ldots$ and $F\in\calS$.
Indeed, for any $n\in\{2,3,\ldots,\infty\}$,
\begin{eqnarray*}
\P(\xi_1+\xi_2\in x+n\Delta)
&=& \sum_{k=0}^{n-1}\P(\xi_1+\xi_2\in x+kT+\Delta)\\
&\sim& 2\sum_{k=0}^{n-1} \P(\xi\in x+kT+\Delta)
= 2\P(\xi\in x+n\Delta).
\end{eqnarray*}

\remark In [\ref{CNW}], the authors consider the class
of distributions concentrated on the integers and such that
$F(\{n+1\})\sim F(\{n\})$ and $F^{*2}(\{n\})\sim 2F(\{n\})$
as $n\to\infty$.
These distributions are $\Delta$-subexponential
with $\Delta=(0,1]$.

\begin{Proposition}\label{inter.integral}
Assume $F[0,\infty )=1$ and $F\in\calL_\Delta$ for some $\Delta$.
Let $\xi_1$ and $\xi_2$ be two i.i.d. random variables with
distribution $F$.
The following assertions are equivalent:

{\rm(i)} $F\in\calS_\Delta$;

{\rm(ii)} there exists a function $h$ such that $h(x)\to\infty$,
$h(x)<x/2$, and $F(x-y+\Delta)\sim F(x+\Delta)$
as $x\to\infty$ uniformly in $|y|\le h(x)$,
\begin{eqnarray}\label{o.small}
{\bf P}(\xi_1+\xi_2\in x{+}\Delta,\ \xi_1>h(x),\ \xi_2>h(x))
&=& o(F(x{+}\Delta))\mbox{ as }x\to\infty;
\end{eqnarray}

{\rm(iii)} the relation {\rm(\ref{o.small})} holds
for every function $h$ such that $h(x)\to\infty$.
\end{Proposition}

\proof. Note that if (\ref{o.small}) is valid for some
$h(x)$, then it follows for any $h_1\ge h$.
For $h(x)<x/2$, the probability of the event
$B=\{\xi_1+\xi_2\in x+\Delta\}$ is equal to
$$
\P(B,\ \xi_1\le h(x))+\P(B,\ \xi_2\le h(x))
+\P(B,\ \xi_1>h(x),\ \xi_2>h(x))
$$
and the conclusions of the proposition follow from
\begin{eqnarray*}
\lefteqn{\P (B,\ \xi_1\le h(x)) = \P (B,\ \xi_2\le h(x))}\\
&&\hspace{20mm} =\int_{0}^{h(x)} F(x{-}y{+}\Delta)F(dy)
\sim F(x{+}\Delta)\int_{0}^{h(x)} F(dy) \sim F(x{+}\Delta).
\end{eqnarray*}

Now we prove that the class $\calS_\Delta$
is closed under a certain local tail equivalence relation.

\begin{Lemma}\label{closure}
Assume that $F\in\calS_\Delta$ for some $\Delta$.
If the distribution $G$ on ${\bf R}^+$ belongs to
$\calL_\Delta$ and
\begin{eqnarray}\label{lim.inf.sup}
0<\liminf_{x\to\infty} \frac{G(x+\Delta)}{F(x+\Delta)}
&\le& \limsup_{x\to\infty} \frac{G(x+\Delta)}{F(x+\Delta)}
<\infty,
\end{eqnarray}
then $G\in\calS_\Delta$. In particular, $G\in\calS_\Delta$,
provided $G(x+\Delta)\sim cF(x+\Delta)$ as $x\to\infty$ for
some $c\in(0,\infty)$.
\end{Lemma}

\proof. Take a function $h(x)\to\infty$ such that
$h(x)<x/2$ and $G(x-y+\Delta)\sim G(x+\Delta)$
as $x\to\infty$ uniformly in $|y|\le h(x)$.
Let $\zeta_1$ and $\zeta_2$ be independent random
variables with common distribution $G$. By Proposition
\ref{inter.integral}(ii), it is sufficient to prove that
\begin{eqnarray*}
I &\equiv&
\P(\zeta_1+\zeta_2\in x+\Delta,\ \zeta_1>h(x),\ \zeta_2>h(x))
= o(G(x+\Delta)).
\end{eqnarray*}
We have
\begin{eqnarray*}
I &=& \int_{h(x)}^{x-h(x)} G(x-y+\Delta) G(dy)
+ \int_{x-h(x)}^{x-h(x)+T}
{\bf P}(\zeta_1\in x-y+\Delta,\ \zeta_1>h(x)) G(dy)\\
&\equiv& I_{1}+I_{2},
\end{eqnarray*}
where
\begin{eqnarray*}
I_{2} &\le& {\bf P}(\zeta_1>h(x))
G(x-h(x)+\Delta) = o(G(x+\Delta))
\end{eqnarray*}
and, by condition (\ref{lim.inf.sup}), for some
$c_1<\infty$ and for all sufficiently large $x$,
\begin{eqnarray*}
I_{1} &\le& c_1 \int_{h(x)}^{x-h(x)} F(x-y+\Delta) G(dy)\\
&\le& c_1\P (\zeta_1+\xi_2\in x+\Delta,\ \zeta_1>h(x),\ \xi_2>h(x))\\
&=& c_1\int_{h(x)}^{x-h(x)} G(x-y+\Delta)F(dy) +
c_1\int_{x-h(x)}^{x-h(x)+T}
{\bf P}(\zeta_1\in x-y+\Delta,\zeta_1>h(x))F(dy)
\end{eqnarray*}
Here $\xi_1$ and $\xi_2$ are independent random
variables with common distribution $F$.
Hence, by using the same arguments as before
and Proposition \ref{inter.integral}(iii),
\begin{eqnarray*}
I_{1} &\le& c_1^2\int_{h(x)}^{x-h(x)} F(x-y+\Delta)F(dy) +
c_1\P (\zeta_1 > h(x)) F(x-h(x)+\Delta)\\
&\le& c_1^2 {\bf P}(\xi_1+\xi_2\in x+\Delta,\ \xi_1\ge h(x),\ \xi_2\ge h(x))
+ o(F(x+\Delta))\\
&=& o(F(x+\Delta))=o(G(x+\Delta)).
\end{eqnarray*}

\begin{Proposition} \label{F.1.F.2}
Assume that $F\in\calS_\Delta$ for some $\Delta$.
Let $G_1$, $G_2$ be two distributions on ${\bf R}_+$
such that $G_1(x+\Delta)/F(x+\Delta)\to c_1$ and
$G_2(x+\Delta)/F(x+\Delta)\to c_2$ as $x\to\infty$,
for some constants $c_1$, $c_2\ge0$. Then
$$
\frac{(G_1*G_2)(x+\Delta)}{F(x+\Delta)} \ \to\ c_1+c_2
\quad\mbox{as } x\to\infty.
$$
If $c_1+c_2>0$ then, by Lemma {\rm\ref{closure}},
$G_1*G_2\in\calS_\Delta$.
\end{Proposition}

\proof. Take two independent random variables $\zeta_1$ and
$\zeta_2$ with distributions $G_1$ and $G_2$.
Take a function $h$ as before. The probability of the event
$B=\{\zeta_1+\zeta_2 \in x+\Delta\}$ is equal to the sum
$$
\P(B,\ \zeta_1\le h(x))+\P(B,\ \zeta_2\le h(x))+
\P(B,\ \zeta_1>h(x),\ \zeta_2>h(x)).
$$
We have that (see the proof of Proposition
\ref{long.tailed}), as $x\to\infty$,
$$
\frac{\P(B,\ \zeta_1\le h(x))}{F(x+\Delta)} \to c_2,\qquad
\frac{\P(B,\ \zeta_2\le h(x))}{F(x+\Delta)} \to c_1.
$$
Following the arguments of Lemma \ref{closure},
we obtain that
$$
\P (B,\ \zeta_1>h(x),\ \zeta_2>h(x)) = o(F(x+\Delta)).
$$
The proposition is proved.

By induction, Proposition \ref{F.1.F.2} implies the following

\begin{Corollary}\label{n.ge.2}
Assume that $F\in\calS_\Delta$ for some $T\in(0,\infty]$
and $G(x+\Delta)/F(x+\Delta)\to c\ge0$ as $x\to\infty$.
Then for any $n\ge2$,
$G^{*n}(x+\Delta)/F(x+\Delta) \to nc$ as $x\to\infty$.
If $c>0$, then $G^{*n}\in\calS_\Delta$.
\end{Corollary}

Let $\{\xi_n\}$ and $\{\zeta_n\}$ be two sequences of
i.i.d.\ non-negative random variables with common
distributions $F(B)=\P(\xi_1\in B)$ and
$G(B)=\P(\zeta_1\in B)$ respectively.
Put $S_n=\zeta_1+\cdots+\zeta_n$.

\begin{Proposition} \label{majorant}
Assume that $F\in\calS_\Delta$ for some $\Delta$
and $G(x+\Delta)=O(F(x+\Delta))$ as $x\to\infty$.
Then, for any $\varepsilon>0$, there exist
$x_0=x_0(\varepsilon)>0$ and $V(\varepsilon)>0$ such that,
for any $x>x_0$ and for any $n\ge 1$,
\begin{eqnarray*}
G^{*n}(x+\Delta) &\le&
V(\varepsilon)(1+\varepsilon)^n F(x+\Delta).
\end{eqnarray*}
\end{Proposition}

\proof. For $x_0\ge 0$ and $k\ge 1$, put
$$
A_k\equiv A_k(x_0) = \sup_{x>x_0}
\frac{G^{*k}(x+\Delta)}{F(x+\Delta)}.
$$
Take any $\varepsilon>0$. Following the arguments of Lemma
\ref{closure}, we conclude the relation, as $x\to\infty$,
$$
{\bf P}(\xi_1+\zeta_2\in x+\Delta,\ \xi_1>h(x),\ \zeta_2>h(x))
=o(F(x+\Delta)).
$$
Hence, there exists $x_0$ such that, for any $x>x_0$,
$$
{\bf P}(\xi_1+\zeta_2\in x+\Delta,\ \zeta_2\le x-x_0)
\le(1+\varepsilon/2)F(x+\Delta).
$$

For any $n>1$ and $x>x_0$,
\begin{eqnarray*}
{\bf P}(S_n\in x+\Delta) &=&
{\bf P} (S_n\in x+\Delta,\ \zeta_n \le x-x_0)\\
&& \hspace{10mm} + {\bf P}(S_n\in x+\Delta,\ \zeta_n>x-x_0)
\equiv P_1(x)+P_2(x),
\end{eqnarray*}
where, by the definition of $A_{n-1}$ and $x_0$,
\begin{eqnarray}\label{P.1.eps}
P_1(x) &=& \int_0^{x-x_0}
{\bf P}(S_{n-1}\in x-y+\Delta) {\bf P}(\zeta_n\in dy)\nonumber\\
&\le& A_{n-1}\int_0^{x-x_0} F(x-y+\Delta){\bf P}(\zeta_n\in dy)\nonumber\\
&=& A_{n-1} {\bf P}(\xi_1+\zeta_n\in x+\Delta,\
\zeta_n\le x-x_0)
\le A_{n-1}(1+\varepsilon/2)F(x+\Delta).
\end{eqnarray}
Further,
\begin{eqnarray*}
P_2(x) &=& \int_0^{x_0+T}
{\bf P}(\zeta_n\in x-y+\Delta,\ \zeta_n>x-x_0){\bf P}(S_{n-1}\in dy)\\
&\le& \sup_{0<t\le x_0}{\bf P}(\zeta_n\in x-t+\Delta)
\int_0^{x_0+T} {\bf P}(S_{n-1}\in dy)
\le \sup_{0<t\le x_0}{\bf P}(\zeta_n\in x-t+\Delta).
\end{eqnarray*}
Thus, if $x>2x_0$, then
\begin{eqnarray*}
P_2(x) &\le& A_1\sup_{0<t\le x_0}F(x-t+\Delta)
\le A_1L_1F(x+\Delta),
\end{eqnarray*}
where
\begin{eqnarray*}
L_1 &=& \sup_{0<t\le x_0,\ y>2x_0}
\frac{F(y-t+\Delta)}{F(y+\Delta)}.
\end{eqnarray*}
If $x_0<x\le 2x_0$, then $P_2(x)\le1$ implies
\begin{eqnarray*}
\frac{P_2(x)}{F(x+\Delta)} &\le&
\frac{1}{\inf_{x_0<x\le 2x_0}F(x+\Delta)}
\equiv L_2.
\end{eqnarray*}
Since $F\in\calL_\Delta$, both $L_1$ and $L_2$ are finite
for $x_0$ sufficiently large. Put $R=A_1L_1+L_2$.
Then, for any $x>x_0$,
\begin{eqnarray} \label{P.2.eps}
P_2(x) &\le& RF(x+\Delta).
\end{eqnarray}
It follows from (\ref{P.1.eps}) and (\ref{P.2.eps})
$A_n \le A_{n-1}(1+\varepsilon/2)+R$ for $n>1$.
Therefore, an induction argument yields:
\begin{eqnarray*}
A_n &\le& A_1(1+\varepsilon/2)^{n-1} +
R\sum_{l=0}^{n-2} (1+\varepsilon/2)^l
\le Rn(1+\varepsilon/2)^{n-1}.
\end{eqnarray*}
This implies the conclusion of the proposition.

Let us consider now some random time $\tau$
with distribution $p_n={\bf P}(\tau=n)$, $n\ge0$
which is independent of $\{\zeta_n\}$.
Then the distribution of the randomly
stopped sum $S_\tau$ is equal to
$$
{\bf P}(S_\tau\in B)=\sum_{n\ge0} p_n G^{*n}(B).
$$

\begin{Theorem}\label{series}
Let $0<T\le\infty$. Assume $F[0,\infty )=1$,
$G(x+\Delta)/F(x+\Delta)\to c\ge0$ as $x\to\infty$,
and ${\bf E}\tau<\infty$.

{\rm(i)} If $F\in\calS_\Delta$ and
${\bf E}(1+\delta)^\tau<\infty$
for some $\delta>0$, then
\begin{eqnarray}\label{series.equiv}
\frac{{\bf P}(S_\tau\in x+\Delta)}{F(x+\Delta)}
&\to& c\cdot{\bf E}\tau\quad\mbox{as }x\to\infty.
\end{eqnarray}

{\rm(ii)} If {\rm(\ref{series.equiv})} holds,
$c>0$, $p_n>0$ for some $n\ge2$, and, in the case of
a finite $\Delta$, $F\in\calL_\Delta$, then $F\in\calS_\Delta$.
\end{Theorem}

\proof\ of (i) follows from Corollary \ref{n.ge.2}, Proposition
\ref{majorant}, and the dominated convergence theorem.

We prove the second assertion.
First, for any $n\ge2$,
\begin{eqnarray}\label{nth.conv.below}
\liminf_{x\to\infty}\frac{G^{*n}(x+\Delta)}{G(x+\Delta)}
&\ge& n.
\end{eqnarray}
Indeed, if $\Delta=(0,\infty)$, then (\ref{nth.conv.below})
follows from Lemma 1 in [\ref{Ch}].
If the interval $\Delta$ is finite, $F\in\calL_\Delta$,
and $c>0$, then $G\in\calL_\Delta$
and (\ref{nth.conv.below})
follows from Corollary \ref{long.tailed.n}.

If $p_n>0$ for some $n\ge2$, then it follows from
(\ref{nth.conv.below}) and (\ref{series.equiv}) that
\begin{eqnarray}\label{some.n}
G^{*n}(x+\Delta) &\sim& nG(x+\Delta)
\quad\mbox{as }x\to\infty
\end{eqnarray}
(the proof is a straightforward argument by contradiction).

If $\Delta=(0,\infty)$, then (\ref{some.n}) implies
the subexponentiality of $G$, by Lemma 7 in [\ref{EGV}].
If $\Delta$ is a finite interval and $F\in\calL_\Delta$,
then $G\in\calL_\Delta$ and, by Corollary \ref{long.tailed.n}, the convolution
$G^{*(n-1)}$ belongs to the class $\calL_\Delta$ too.
Thus, by Proposition \ref{long.tailed},
\begin{eqnarray*}
n=\limsup_{x\to\infty}\frac{G^{*n}(x+\Delta)}{G(x+\Delta)}
&=& \limsup_{x\to\infty}\frac{(G*G^{*(n-1)})(x+\Delta)}{G(x+\Delta)}\\
&\ge&
1+\limsup_{x\to\infty}\frac{G^{*(n-1)}(x+\Delta)}{G(x+\Delta)}.
\end{eqnarray*}
By induction we deduce from this estimate that
\begin{eqnarray*}
\limsup_{x\to\infty}\frac{G^{*2}(x+\Delta)}{G(x+\Delta)}
&\le& 2,
\end{eqnarray*}
which implies the $\Delta$-subexponentiality of $G$.
Now $F\in\calS_{\Delta}$ by Lemma \ref{closure}.

In Theorem \ref{series}, assertion (i) is valid for any
$\Delta$-subexponential distribution. For a fixed
distribution $F$, the condition
${\bf E}(1+\delta)^\tau<\infty$ may be substantially
weakened. We can illustrate that by the following example.
Consider the case of the infinite interval $\Delta=(0,\infty)$.
Assume that $G=F$ and there exist finite positive constants $c$ and
$\alpha$ such that $\overline F(x/n)\le cn^\alpha\overline F(x)$
for any $x>0$ and $n\ge1$ (for instance, the Pareto
distribution with parameter $\alpha$ satisfies this condition).
Then ${\bf P}(S_\tau>x)\sim {\bf E}\tau\cdot \overline F(x)$
as $x\to\infty$ provided ${\bf E}\tau^{1+\alpha}$ is finite,
as follows by combining dominated convergence with
\begin{eqnarray*}
{\bf P}(S_n>x) &\le& {\bf P}(n\cdot\max_{k\le n}\zeta_k>x)
\le n{\bf P}(\zeta_1>x/n) \le n^{1+\alpha}\overline F(x).
\end{eqnarray*}

Proposition \ref{majorant} implies also the following corollary.
For $x\ge0$, put $\eta(x)=\min\{n\ge1:S_n>x\}$ and
$\chi(x)=S_{\eta(x)}-x$. As earlier,
let $\tau$ be a non-negative integer-valued random variable
which does not depend on $\zeta$'s.

\begin{Corollary} \label{overshoot}
Assume that $G\in\calS_\Delta$ and
${\bf E}(1+\delta)^\tau<\infty$ for some $\delta>0$.
Then ${\bf P}(\chi(x)\in y+\Delta,\ \eta(x)\le\tau)
\sim {\bf E}\tau\cdot G(x+y+\Delta)$ as $\min(x,y)\to\infty$.
\end{Corollary}

\proof. Let $h$ be such that $h(y)\le y/2$,
$h(y)\uparrow\infty$ as $y\to\infty$, and
$G(y+t+\Delta)\sim G(y+\Delta)$ uniformly in $|t|\le h(y)$.
Put $z=\min(h(y),x)$. For any $n\ge2$,
\begin{eqnarray*}
{\bf P}(\chi(x)\in y+\Delta,\ \eta(x)=n)
&=& {\bf P}(S_{n-1}\le x,\ S_n\in x+y+\Delta)\\
&\ge& \int_0^z {\bf P}(S_{n-1}\in dt)G(x+y-t+\Delta)
\sim G(x+y+\Delta).
\end{eqnarray*}
On the other hand,
\begin{eqnarray*}
\lefteqn{{\bf P}(S_{n-1}\le x,\ S_n\in x+y+\Delta)}\\
&&\hspace{10mm} = {\bf P}(S_n\in x+y+\Delta)
-{\bf P}(S_{n-1}>x,\ S_n\in x+y+\Delta)\\
&&\hspace{10mm} \le {\bf P}(S_n\in x+y+\Delta)
-{\bf P}(\zeta_n\le z,\ S_n\in x+y+\Delta)\\
&&\hspace{10mm} = (n+o(1))G(x+y+\Delta)-(n-1+o(1))G(x+y+\Delta).
\end{eqnarray*}
Thus, for any fixed $n\ge1$,
$$
{\bf P}(\chi(x)\in y+\Delta,\ \eta(x)=n)\sim G(x+y+\Delta).
$$
Now Proposition \ref{majorant} and the dominated convergence theorem
complete the proof, since ${\bf P}(\chi(x)\in y+\Delta,\
\eta(x)=n)\le {\bf P}(S_n\in x+y+\Delta)$ and
\begin{eqnarray*}
{\bf P}(\chi (x)\in y+\Delta,\ \eta(x)\le\tau)
&=& \sum_{k=1}^\infty {\bf P}(\tau=k)\sum_{n=1}^k
{\bf P}(\chi(x)\in y+\Delta,\ \eta(x)=n)\\
&\sim& G(x+y+\Delta) \sum_{k=1}^\infty k{\bf P}(\tau=k).
\end{eqnarray*}

\mysection{Distributions with subexponential densities}\label{Local}

Similar results (with similar proofs!) hold
for densities of absolutely continuous distributions.
More precisely, in this Section we consider a class
of distributions $\{F\}$ with the following property:
each distribution $F$ has a density $f(x)$ for all
sufficiently large values of $x$, i.e., for a certain
$\widehat x=\widehat x(F)$ and for any Borel set
$B\subseteq[\widehat x,\infty)$,
\begin{eqnarray*}
F(B) &=& \int_B f(y)dy.
\end{eqnarray*}

We say that a density $f$ on $[\widehat x(F),\infty)$ is
{\it long-tailed} (and write $f\in\calL$) if the function
$f(x)$ is bounded on $[\widehat x,\infty)$,
$f(x)>0$ for all sufficiently large $x$, and
$f(x+t)\sim f(x)$ as $x\to\infty$ uniformly in $t\in[0,1]$.
In particular, if $f\in \calL$, then $f(x)\to 0$ as $x\to\infty$.

A distribution $F$ on ${\bf R}_+$ with a density $f(x)$
on $[\widehat x,\infty)$ is said to belong
to the class $\calS_{ac}$
(the density $f$ is {\it subexponential})
if $f\in\calL$ and, as $x\to\infty$,
\begin{eqnarray*}
f^{*2}(x) \equiv 2\int_0^{\widehat x} f(x-y)F(dy)
+\int_{\widehat x}^{x-\widehat x} f(x-y)f(y)dy
&\sim& 2f(x).
\end{eqnarray*}

Typical examples of $\calS_{ac}$ are given by the
Pareto, lognormal, and Weibull (with parameter between 0 and 1)
distributions (for the proof, see Section \ref{Sufficient}).
Note that distribution with subexponential density is
$\Delta$-subexponential for any $0<T\le\infty$.

\begin{Proposition}\label{long.tailed.dens}
Let $F$ and $G$ have  densities $f$ and $g$
on $[\widehat x,\infty)$ belonging to the class $\calL$.
Then the density $f*g$ of the convolution $F*G$ is long-tailed and
\begin{eqnarray}\label{long.tailed.below.dens}
\liminf_{x\to\infty} \frac{(f*g)(x)}{f(x)+g(x)} &\ge& 1.
\end{eqnarray}
In particular, if $f\in\calL$, then $f^{*n}\in\calL$ and
$\liminf_{x\to\infty}f^{*n}(x)/f(x)\ge n$.
\label{long.tailed.n.dens}
\end{Proposition}

\proof. Take a function $h(x)\uparrow\infty$ such that
$\widehat x\le h(x)<x/2$, $f(x-y)\sim f(x)$ and $g(x-y)\sim g(x)$
as $x\to\infty$ uniformly in $|y|\le h(x)$. Then
\begin{eqnarray*}
(f*g)(x+t) &=& \int_{-\infty}^{x-h(x)}f(x+t-y)G(dy)
+\int_{-\infty}^{h(x)} g(x+t-y)F(dy)\\
&& +\int_{x-h(x)}^{x+t-h(x)} f(x+t-y)g(y)dy
\equiv I_1(x,t)+I_2(x,t)+I_3(x,t).
\end{eqnarray*}
Now the conclusion of the proposition follows from
$I_1(x,t)\sim I_1(x,0)$ and $I_2(x,t)\sim I_2(x,0)$
as $x\to\infty$ uniformly in $t\in(0,1]$ and the estimate
\begin{eqnarray*}
I_3(x,t) &\le& \sup_{y\in[h(x),h(x)+t)}f(y)
\int_{x-h(x)}^{x+t-h(x)} g(y)dy
\sim tf(h(x))g(x)=o(g(x)).
\end{eqnarray*}

\begin{Proposition}\label{inter.integral.dens}
Assume that the distribution $F$ on ${\bf R}^+$ has a
density $f\in\calL$ on $[\widehat x,\infty)$.
Then the following assertions are equivalent:

{\rm(i)} the density $f$ is subexponential;

{\rm(ii)} for some function $h$ such that $h(x)\to\infty$,
$h(x)<x/2$, and $f(x-y)\sim f(x)$
as $x\to\infty$ uniformly in $|y|\le h(x)$,
\begin{eqnarray}\label{o.small.dens}
\int_{h(x)}^{x-h(x)} f(x-y)f(y) dy
&=& o(f(x))\mbox{ as }x\to\infty;
\end{eqnarray}

{\rm(iii)} the relation {\rm(\ref{o.small.dens})} holds
for every function $h$ such that $h(x)\to\infty$.
\end{Proposition}

\proof. For $\widehat x\le h(x)<x/2$,
$$
f^{*2}(x) = 2\int_{0}^{h(x)} f(x-y)F(dy)
+\int_{h(x)}^{x-h(x)} f(x-y)f(y) dy.
$$
Here the first integral is equivalent to $f(x)$
as $x\to\infty$. This completes the proof.

\begin{Lemma}\label{closure.dens}
Let $f$ be a subexponential density on $[\widehat x,\infty)$.
Assume that the density $g$ on $[\widehat x,\infty)$
is long-tailed and
\begin{eqnarray*}
0<\liminf_{x\to\infty}g(x)/f(x)
&\le& \limsup_{x\to\infty}g(x)/f(x) <\infty.
\end{eqnarray*}
Then $g$ is subexponential too. In particular,
$g\in\calS_{ac}$, given $g(x)\sim cf(x)$
as $x\to\infty$ for some $c\in(0,\infty)$.
\end{Lemma}

\proof. The result follows by Proposition
\ref{inter.integral.dens}(iii). Indeed,
one can choose $c_1<\infty$ such that
$g(x)\le c_1f(x)$ for all sufficiently large $x$ and
\begin{eqnarray*}
\int_{h(x)}^{x-h(x)} g(x-y)g(y)dy
&\le& c_1^2\int_{h(x)}^{x-h(x)} f(x-y)f(y)dy.
\end{eqnarray*}

\begin{Proposition} \label{F.1.F.2.dens}
Let $f$ be a subexponential density on $[\widehat x,\infty)$.
Let $f_1$, $f_2$ be two densities on $[\widehat x,\infty)$
such that $f_1(x)/f(x)\to c_1$ and $f_2(x)/f(x)\to c_2$
as $x\to\infty$, for some constants $c_1$, $c_2\ge0$. Then
$$
\frac{(f_1*f_2)(x)}{f(x)} \ \to\ c_1+c_2
\quad\mbox{as } x\to\infty.
$$
If $c_1+c_2>0$ then, by Lemma {\rm\ref{closure.dens}},
the convolution $f_1*f_2$ is a subexponential density.
\end{Proposition}

\proof. Take a function $h$ as before. Then
\begin{eqnarray*}
f_1*f_2(x) &=& \int_{0}^{h(x)} f_1(x-y)F_2(dy)
+\int_{0}^{h(x)} f_2(x-y)F_1(dy)
+\int_{h(x)}^{x-h(x)} f_1(x-y)f_2(y) dy\\
&\equiv& I_1(x)+I_2(x)+I_3(x).
\end{eqnarray*}
We have $I_1(x)/f(x)\to c_1$ and
$I_2(x)/f(x)\to c_2$ as $x\to\infty$.
Finally,
$$
I_3(x) \le
(c_1c_2+o(1))\int_{h(x)}^{x-h(x)} f(x-y)f(y) dy = o(f(x)),
$$
which completes the proof.

\begin{Corollary}\label{n.ge.2.dens}
Assume that $F\in\calS_{ac}$.
Then for any $n\ge2$, $f^{*n}(x)\sim nf(x)$ as
$x\to\infty$ and $F^{*n}\in\calS_{ac}$.
\end{Corollary}

Let $\{\xi_n\}$ be a sequence of i.i.d.\ non-negative
random variables with a common distribution
$F(B)=\P(\xi_1\in B)$. Put $S_n=\xi_1+\cdots+\xi_n$.

\begin{Proposition}\label{majorant.dens}
Assume that $F\in\calS_{ac}$.
Then, for any $\varepsilon>0$, there exist
$x_0=x_0(\varepsilon)\ge\widehat x$ and $V(\varepsilon)>0$ such that,
for any $x>x_0$ and for any integer $n\ge 1$,
\begin{eqnarray*}
f^{*n}(x) &\le& V(\varepsilon)(1+\varepsilon)^n f(x).
\end{eqnarray*}
\end{Proposition}

\proof. For $x_0\in{\bf R}_+$ and $k\ge 1$, put
$A_k\equiv A_k(x_0) = \sup_{x>x_0} f^{*k}(x)/f(x)$.
Take any $\varepsilon >0$. Fix an integer
$j$ such that $(j+1+\varepsilon)^{1/j}<1+\varepsilon/2$.
By Corollary \ref{n.ge.2.dens},
$A_{j+1}(x_0) \to j+1$ as $x_0\to\infty$.
Choose $x_0\ge\widehat x$ such that
$$
A_{j+1}(x_0)\le (j+1+\varepsilon)<(1+\varepsilon/2)^j.
$$

For any $n>j$ and $x>2x_0$,
\begin{eqnarray*}
f^{*n}(x) &=&
\int_{0}^{x-x_0} f^{*(n-j)}(x-y)F^{*j}(dy)
+\int_{0}^{x_0} f^{*j}(x-y)F^{*(n-j)}(dy)
\equiv I_1(x)+I_2(x),
\end{eqnarray*}
where, by the definition of $A_{n-j}$ and $A_{j+1}$,
\begin{eqnarray}\label{I.1.eps.dens}
I_1(x) &\le&
A_{n-j}\int_{0}^{x-x_0} f(x-y)F^{*j}(dy)
\le A_{n-j} f^{*(j+1)}(x) \le A_{n-j} A_{j+1}f(x)
\end{eqnarray}
and
\begin{eqnarray}\label{I.2.eps.dens}
I_2(x) &\le& \max_{0<t\le x_0}f^{*j}(x-t)
\le A_j \max_{0<t\le x_0}f(x-t)\le A_jL_1f(x),
\end{eqnarray}
where $L_1=\sup_{0<t\le x_0, y>2x_0} f(y-t)/f(y)$.
If $x_0<x\le 2x_0$, then
\begin{eqnarray}\label{le.2x_0}
\frac{f^{*n}(x)}{f(x)} &\le&
\frac{\sup_{x\in(x_0,2x_0]}f^{*n}(x)}
{\inf_{x_0<x\le 2x_0}f(x)}
\equiv L_2.
\end{eqnarray}
Since $f\in\calL$, we may choose $x_0$ such that $L_1$ and
$L_2<\infty$. Put $R=A_jL_1+L_2$. It follows from
(\ref{I.1.eps.dens})--(\ref{le.2x_0}) that, for any $x>x_0$,
\begin{eqnarray*}
f^{*n}(x) &\le& (A_{n-j}A_{j+1}+R)f(x).
\end{eqnarray*}
Hence, for $n>j$, $A_n \le A_{n-j}A_{j+1} +R$.
The remaining part of the proof is the same as that of
Proposition \ref{majorant}.

\begin{Theorem}\label{series.dens}
Let $\{p_n\}_{n\ge1}$ be a non-negative sequence such that
$m_p\equiv\sum_{n\ge1}np_n$ is finite. Denote
$$
g(x)=\sum_{n\ge1} p_n f^{*n}(x).
$$

{\rm(i)} If a distribution $F$ has a subexponential
density $f$ on $[\widehat x,\infty)$ and
$\sum_{n\ge1}(1+\delta)^n p_n<\infty$ for some $\delta>0$, then
\begin{eqnarray}\label{series.equiv.dens}
g(x) &\sim& m_pf(x) \quad\mbox{as }x\to\infty.
\end{eqnarray}

{\rm(ii)} If equivalence {\rm(\ref{series.equiv.dens})} holds,
$p_1<1$, $F[0,\infty )=1$, and $f\in\calL$, then $F\in\calS_{ac}$.
\end{Theorem}

\proof. Assertion (i) follows from Corollary
\ref{n.ge.2.dens}, Proposition \ref{majorant.dens},
and the dominated convergence theorem.
We prove the second assertion.
By Lemma \ref{long.tailed.n.dens}, for any $n\ge2$,
\begin{eqnarray*}
\liminf_{x\to\infty} f^{*n}(x)/f(x) &\ge& n.
\end{eqnarray*}
If $p_n>0$ for some $n\ge2$, then this estimate
and (\ref{series.equiv.dens}) imply that
\begin{eqnarray}\label{some.n.dens}
f^{*n}(x) &\sim& nf(x) \quad\mbox{as }x\to\infty.
\end{eqnarray}

By Proposition \ref{long.tailed.n.dens}, $f^{*(n-1)}\in\calL$ and
\begin{eqnarray*}
n=\limsup_{x\to\infty}\frac{f^{*n}(x)}{f(x)}
&=& \limsup_{x\to\infty}\frac{(f*f^{*(n-1)})(x)}{f(x)}
\ge 1+\limsup_{x\to\infty}\frac{f^{*(n-1)}(x)}{f(x)}.
\end{eqnarray*}
By induction we deduce from this estimate that
$\limsup_{x\to\infty} f^{*2}(x)/f(x)\le 2$,
which implies the subexponentiality of $f$.

\mysection{Sufficient conditions for
$\Delta$-subexponentiality\\
and subexponentiality of densities}
\label{Sufficient}

The sufficient conditions for distributions to be
subexponential are well-known (see, e.g., [\ref{T},
\ref{RSST} Section 2.5.3]).
In this Section, we propose similar conditions for
distributions to belong either to $\calS_\Delta$ for
a finite $T$, or to $\calS_{ac}$.

\begin{Proposition}\label{suff.for.Delta}
Let a distribution $F$ on ${\bf R}^+$ belong to the class
$\calL_\Delta$ for some finite $T>0$. Assume that
there exist $c>0$ and $x_0<\infty$ such that
$F(x+t+\Delta)\ge cF(x+\Delta)$ for any
$t\in(0,x]$ and $x>x_0$. Then $F\in\calS_\Delta$.
\end{Proposition}

\proof. Let a function $h(x)\to\infty$ be such that
$h(x)<x/2$. Then
\begin{eqnarray*}
{\bf P}(\xi_1+\xi_2\in
x{+}\Delta,\xi_1>h(x),\xi_2>h(x))
&\le& 2\int_{h(x)}^{x/2+T} F(x-y+\Delta)F(dy)\\
&& \hspace{-60mm}\le 2(c+o(1))\int_{h(x)}^{x/2+T} F(x+\Delta)F(dy)
=o(F(x+\Delta))
\end{eqnarray*}
as $x\to\infty$. Applying now Lemma \ref{inter.integral}(ii)
we conclude that $F\in\calS_\Delta$.

The Pareto distribution (with the tail
$\overline F(x)=x^{-\alpha}$, $\alpha>0$, $x\ge1$)
satisfies conditions of Proposition \ref{suff.for.Delta}.
The same is true for any distribution $F$ such that
${\bf P}(\xi\in x+\Delta)$ is regularly
varying at infinity,
i.e., for $F(x+\Delta)\sim x^{-\alpha}l(x)$,
where $l(x)$ is slowly varying at infinity.

\begin{Proposition}\label{suff.for.Weibull-type}
Let a distribution $F$ on ${\bf R}^+$ belong to the class
$\calL_\Delta$ for some finite $\Delta$.
Let there exist $x_0$ such that the function
$g(x)\equiv -\ln F(x+\Delta)$ is concave for $x\ge x_0$.
Let there exist a function $h(x)\uparrow\infty$ as $x\to\infty$
such that $F(x+t+\Delta)\sim F(x+\Delta)$ uniformly in
$|t|\le h(x)$ and $x F(h(x)+\Delta )\to 0$.
Then $F\in\calS_\Delta$.
\end{Proposition}

\proof. Due to Lemma \ref{closure},
without loss of generality assume $x_0=0$.
Since $g(x)$ is concave, the minimum of the sum
$g(x-y)+g(y)$ on the interval $y\in[h(x),x-h(x)]$
is equal to $g(x-h(x))+g(h(x))$. Therefore,
\begin{eqnarray*}
\int_{h(x)}^{x-h(x)} F(x-y+\Delta)F(dy) &\le&
c_1\int_{h(x)}^{x-h(x)} F(x-y+\Delta)F(y+\Delta)dy\\
&=& c_1\int_{h(x)}^{x-h(x)} e^{-(g(x-y)+g(y))}dy
\le c_1xe^{-(g(x-h(x))+g(h(x)))}
\end{eqnarray*}
Since $e^{-g(x-h(x))}\sim e^{-g(x)}$,
\begin{eqnarray*}
\int_{h(x)}^{x-h(x)} F(x-y+\Delta)F(dy)
&=& O(e^{-g(x)}xe^{-g(h(x))}) = o(F(x+\Delta)),
\end{eqnarray*}
which completes the proof.

Consider the Weibull distribution,
$\overline F(x)=e^{-x^\beta}$, $x\ge0$, $\beta\in(0,1)$. Then
$$
F(x+\Delta ) \sim \beta T x^{\beta -1} \exp (-x^{\beta})
\quad \mbox{as} \quad x\to \infty .
$$
Consider the distribution $\widehat F$ with the tail
$\overline{\widehat F}(x)=\min(1,x^{\beta-1}e^{-x^\beta})$.
Let $x_0$ be the unique positive solution to the equation
$x^{1-\beta}=e^{-x^\beta}$. Then the function
$\widehat g(x)=-\ln\widehat F(x+\Delta)$ is concave
for $x\ge x_0$,  and conditions of Proposition
\ref{suff.for.Weibull-type} are satisfied with $h(x)=x^\gamma$,
$\gamma \in (0,1-\beta )$. Therefore,
$\widehat F\in\calS_\Delta$ and, due to Lemma
\ref{closure}, $F\in\calS_\Delta$.

Similarly, one can check that, for the lognormal
distribution with the density
$f(x)=e^{-(\ln x-\ln a)^2/2\sigma^2}/x\sqrt{2\pi\sigma^2}$,
$$
F(x+\Delta ) \sim Tf(x),
$$
the function $g(x)=-\ln(x^{-1}e^{-(\ln x-\ln a)^2/2\sigma^2})
=\ln x+(\ln x-\ln a)^2/2\sigma^2$ is eventually concave,
and conditions of Proposition \ref{suff.for.Weibull-type}
are satisfied with any $h(x)=o(x)$. Thus, $F\in\calS_\Delta$.

Similarly to Propositions \ref{suff.for.Delta} and
\ref{suff.for.Weibull-type} we obtain the following

\begin{Proposition}\label{suff.for.local}
Let a distribution $F$ on ${\bf R}^+$ have a long-tailed
density $f(x)$. Let one of the following conditions hold:

{\rm(i)} there exists $c>0$ such that
$f(y)\ge cf(x)$ for any $y\in(x,2x]$;

{\rm(ii)} the function $g(x)\equiv -\ln f(x)$ is concave
for $x\ge x_0$ and, for some $h(x)\to\infty$,
$f(x+t)\sim f(x)$ uniformly in $|t|\le h(x)$
and $xe^{-g(h(x))}\to 0$.

Then $f$ is subexponential.
\end{Proposition}

The density of the Pareto distribution satisfies
condition (i) of Proposition \ref{suff.for.local}.
The density of the Weibull distribution with parameter
$\beta\in(0,1)$ satisfies condition (ii) of Proposition
\ref{suff.for.local} with $h(x)=\ln^{2/\beta}x$.

{\bf Example 1.} Assume that $\xi$ takes positive integer
values only, $\P(\xi=2k)=\gamma/k^2$ and
$\P(\xi=2k+1)=\gamma/2^k$,
where $\gamma$ is a normalizing constant.
Then $\xi$ has a lattice distribution $F$ with span $1$.
By Proposition \ref{suff.for.Delta}, $F\in\calS_{(0,2]}$,
but $F$ cannot belong to any $\calS_{(0,a]}$
if $a$ is not infinity or an even integer.

{\bf Example 2.} Assume that $\xi$ is a sum of two
independent random variables: $\xi=\eta+\zeta$
where $\eta$ is distributed uniformly on $(-1/8,1/8)$
and $\P(\zeta=k)=\gamma/k^2$. Then the distribution $F$
of $\xi$ is absolutely continuous.
It may be checked that $F\in\calS_{(0,1]}$,
but $F$ cannot belong to any $\calS_{(0,a]}$
if $a$ is not infinity or an integer.

{\bf Example 3.} Consider a long-tailed function $f(x)$
in the range $f(x)\in[1/x^2,2/x^2]$ for any $x > 0$.
Let us choose the function $f$ in such a way that $f$ is
not asymptotically equivalent to a non-increasing function.

For instance, one can define $f$ as follows.
Consider the increasing sequence $x_n=2^{n/4}$.
Put $f(x_{2n})= 1/x_{2n}^2$ and
$f(x_{2n+1})=2/x_{2n+1}^2$. Then assume that $f$
is linear between any two consecutive points.

Consider the lattice distribution $F$ on the set of natural
numbers with $F(\{n\})=f(n)$ for all sufficiently large
integer $n$. Then by Lemma \ref{closure},
$F\in\calS_{(0,1]}$, but $f(n)=F((n-1,n])$
is not asymptotically equivalent
to a non-increasing function.

{\bf Example 4.} Let $G_+$ be the ascending ladder height
distribution of a random walk with increment distribution
$F$. It is shown in [\ref{AKKKT}] that $G_+\in \calS_\Delta$
for all $T<\infty$ when $F$ is non-lattice. However, $G_+$
cannot have a subexponential density when $F$ is singular
(say concentrated on $\{-1,\sqrt{2}\}$) since then also
$G_+$ is singular.

\mysection{Supremum of a random walk}\label{Supremum}

Theorems \ref{series} and \ref{series.dens} give us
a unified approach for obtaining the local and integral
asymptotic theorems for the supremum of a random walk.

Let $\{\xi_n\}$ be a sequence of independent random variables
with a common distribution $F(B)=\P(\xi_n\in B)$ and
$\E\xi_1=-m<0$. Let
$$
F^I(x) \equiv 1-\overline{F^I}(x)
= 1-\min(1,\int_x^{\infty}\overline{F}(y)dy)
$$
denote the integrated-tail distribution function.
It is easy to see that

(a) if $ F$ is long-tailed, then $F^I$ is long-tailed, too;

(b) $F^I$ is long-tailed if and only if
$\overline F(x) =o(\overline{F^I}(x))$ as $x\to\infty$.

Put $S_0=0$, $S_n=\xi_1+\cdots+\xi_n$. By the SLLN,
$M=\sup_{n\ge0} S_n$ is finite with probability 1.
Write $\pi(B) = \P(M\in B)$,
$\pi(x) \equiv \pi (-\infty ,x] = 1-\overline \pi(x)$.

It is well-known (see, e.g. [\ref{A}, \ref{EKM}, \ref{EV}]
and references therein)
that if $F^I \in\calS$, then, as $x\to\infty$,
\begin{equation} \label{integral}
\overline{\pi}(x) \sim \frac{1}{m} \overline{F^I}(x).
\end{equation}
In particular, $\pi\in\calS$.
Korshunov [\ref{K}] proved the converse:
(\ref{integral}) implies $F^I\in\calS$.

Recently, Asmussen et al. [\ref{AKKKT}] proved that if
$F\in\calS\,^*$, i.e. if
$$
\int_0^x \overline{F}(x-y) \overline{F}(y) dy \sim 2 \E \max (\xi_1,0)
\overline{F}(x), \quad x\to\infty ,
$$
then, for any $T\in (0,\infty )$,
\begin{equation} \label{local}
\pi(x+\Delta ) \sim \frac{T}{m}\overline F(x)
\end{equation}
(if the distribution $F$ is lattice then $x$ and $T$
should be restricted to values of the lattice span).
In particular, $\pi\in\calS_\Delta$ for any $0<T<\infty$.

In the lattice case, (\ref{local}) was proved earlier by
Bertoin and Doney [\ref{BD}]. They also sketched a proof
of (\ref{local}) for non-lattice distributions.

It follows from [\ref{FZ}, Theorem 2(b)] that the converse
is also true: if (\ref{local}) holds for any $T\in(0,\infty)$
and $F$ is long tailed, then $F\in\calS\,^*$.

\remark
Since (\ref{local}) holds for any $T>0$,
it implies that, for any $T_0>0$,
\begin{equation} \label{local_uniform}
\pi(x+\Delta ) \sim \frac{1}{m}\int_x^{x+T}\overline F(y)dy
\end{equation}
as $x\to\infty$ uniformly in $T\in[T_0,\infty]$.

One can see that Theorem \ref{series} gives a unified
approach for obtaining (\ref{integral})--(\ref{local}).
We start with the following

\begin{Lemma}\label{first.positive}
Let $v(x)$ be a long-tailed function and let
$$
V(x) \equiv \int_x^\infty v(y)dy.
$$
Assume that $V(0)<\infty$.
For any $n$, define the event
$A_n=\{S_j\le0\mbox{ for all }j\le n\}$ and put $p=\P(M>0)$.
Then, as $x\to\infty$,
\begin{eqnarray*}
\sum_{n=0}^\infty \E(v(x-S_n);\ A_n)
&\sim& \frac{1-p}{m}V(x).
\end{eqnarray*}
\end{Lemma}

\proof. Since $v$ is long-tailed, $V$ is long-tailed, too,
and $v(x) = o(V(x))$.

Assume that the distribution $F$
is non-lattice (the proof in the lattice case is similar).
For $n\ge0$, consider the measures
\begin{eqnarray*}
H_n(B) &=& {\bf P}\{S_j\le0 \mbox{ for any } j\le n,\ S_n\in B\},
\quad B\subseteq(-\infty,0]
\end{eqnarray*}
and the corresponding taboo renewal function
\begin{eqnarray*}
H(B) &=& \sum_{n=0}^\infty H_n(B).
\end{eqnarray*}
It is well-known that, for a non-lattice distribution,
\begin{eqnarray}\label{asympt.of.Psi}
H(y+(0,1]) &\sim& (1-p)/m\quad \mbox{as }y\to-\infty.
\end{eqnarray}
Since
\begin{eqnarray*}
\E(v(x-S_n);\ A_n) &=& \int_{-\infty}^0 v(x-y)H_n(dy)
\end{eqnarray*}
and the function $v(x)$ is long-tailed, we obtain
\begin{eqnarray*}
\sum_{n=0}^\infty \E(v(x-S_n);\ A_n)
&=& \int_0^\infty v(x+y) H(-dy)
\sim \sum_{j=0}^\infty v(x+j) H((-j-1,-j]).
\end{eqnarray*}
Take an integer-valued function $h(x)\to\infty$ such that
$v(x+t)\sim v(x)$ uniformly in $|t|\le h(x)$
and $v(x)h(x)=o(V(x))$. Then, by (\ref{asympt.of.Psi}),
\begin{eqnarray*}
\sum_{n=0}^\infty \E(v(x-S_n);\ A_n)
&\sim& \sum_{j=h(x)}^\infty v(x+j) H((-j-1,-j])\\
&\sim& \frac{1-p}{m}\sum_{j=h(x)}^\infty v(x+j)
\sim \frac{1-p}{m}\int_x^\infty v(y)dy.
\end{eqnarray*}
The proof is complete.

Consider the defective stopping time
$$
\eta = \inf \{ n\ge1 : \ S_n >0 \} \le\infty
$$
and let $\{\psi_n\}$ be i.i.d.\ random variables
with common distribution function
$$
G(x) \equiv  \P(\psi_n\le x)=\P(S_{\eta}\le x~|~\eta<\infty).
$$
It is well-known (see, e.g. Feller [\ref{F}, Chapter 12]) that the
distribution of the maximum $M$ coincides with the distribution
of the randomly stopped sum $\psi_1+\cdots+\psi_\nu$,
where the stopping time $\nu$ is independent of
the sequence $\{\psi_n\}$ and is geometrically
distributed with parameter $p=\P(M>0)<1$, i.e.,
$\P(\nu=k)=(1-p)p^k$ for $k=0$, 1, \ldots. Equivalently,
\begin{eqnarray*}
\P(M\in B) &=& (1-p)\sum_{k=0}^\infty p^k G^{*k}(B).
\end{eqnarray*}
From Borovkov [\ref{B}, Chapter 4, Theorem 10],
if ${F}^I$ is long-tailed, then
\begin{equation} \label{qq}
\overline{G}(x) \sim \frac{1-p}{pm} \overline{F^I}(x).
\end{equation}
For any $T\in(0,\infty)$ and $\Delta=(0,T]$,
if the function $v(x)=F(x+\Delta)$ is long-tailed,
then by Lemma \ref{first.positive},
\begin{eqnarray}\label{qqq}
G(x+\Delta) &=& \P(S_\eta\in x+\Delta)/\P(\eta<\infty)
= \frac{1}{p}\sum_{n=1}^\infty\P(S_n\in x+\Delta,\ \eta=n)\nonumber\\
&\sim& \frac{1-p}{pm} \int_x^\infty F(y+\Delta ) \sim
\frac{(1-p)T}{pm} \overline{F}(x).
\end{eqnarray}
Now (\ref{integral})--(\ref{local}) follow from
(\ref{qq})--(\ref{qqq}) by Theorem \ref{series}.

Similarly, Theorem \ref{series.dens} allows us to get the
asymptotics for the density of $\pi$.

\begin{Theorem}\label{rwdensity}
Assume that $F\in\calS\,^*$ and that the density $f$
on $[x(F),\infty)$ of the distribution $F$ is long-tailed.
Then, as $x\to\infty$, the density
of $\pi$ is equivalent to $\overline{F}(x)/m$.
\end{Theorem}

Indeed, if the density $f$ of the distribution $F$ is long-tailed,
then by Lemma \ref{first.positive} (with $v(x)=f(x)$),
$G$ has a density $g$ on the interval $[\hat{x}(F),\infty )$
which is long-tailed and
\begin{eqnarray*}
g(x) &\sim& \frac{1-p}{pm}\overline F(x).
\end{eqnarray*}
Further, if $F\in\calS\,^*$, then $G$ has a subexponential
density.
The density of  the distribution $\pi$ may be represented as
$$
(1-p)\sum_{k=1}^\infty p^k g^{*k}(x),
$$
and, by Theorem \ref{series.dens}(i), is equivalent to
\begin{eqnarray}\label{rwd}
g(x)(1-p)\sum_{k=1}^\infty kp^k
&\sim& \overline F(x)/m\quad\mbox{as }x\to\infty.
\end{eqnarray}

\remark The result of Theorem \ref{rwdensity} is new. In
[\ref{AKKKT}, Proposition 1], it was claimed that
the same asymptotics may be obtained under different
conditions.

\mysection{The renewal function and the key renewal theorem}
\label{KRT}

Let $G$ be a non-negative measure on $(0,\infty)$.
We will assume throughout that $\theta\le 1$ where
$\theta=G(0,\infty)$. Then the renewal measure
$$
U=\sum_{n=0}^\infty G^{*n}
$$
is well-defined and finite on compact sets.
In addition, if $\theta<1$ then $U$ is a finite measure
(in fact, $U[0,\infty)=$ $(1-\theta)^{-1}$).
See, e.g., [\ref{F}] or [\ref{APQ}]
for this and further basic facts from renewal theory.

Blackwell's renewal theorem states that when $\theta=1$
and $G$ is non-lattice, then $U(x+\Delta)\sim T/\mu_G$
where $\mu_G$ is the mean of $G$.
When   $\theta<1$ and $G$ is light-tailed, it is
easy to see by standard techniques
([\ref{APQ}], VI.5) that $U(x+\Delta)$
decreases exponentially fast.
Callaert \& Cohen [\ref{CC}] gave an asymptotic expression
for a special heavy-tailed case with $\theta<1$, $T=\infty$.
Here is a more complete and
local version. We will say that $G\in\calS_\Delta$ if
$F \in\calS_\Delta$ where $F$ is the probability measure
$G/\theta$.

\begin{Proposition}\label{renfct}
Assume $0<\theta<1$.
If $T<\infty$, assume also that $G\in\calL_\Delta$.
Then $U(x+\Delta)\sim (1-\theta)^{-2}G(x+\Delta)$ as $x\to\infty$
if and only if $G\in \calS_\Delta$.
\end{Proposition}

\proof. By Theorem \ref{series},
\begin{eqnarray*}
U(x+\Delta) &=&
\sum_{n=1}^\infty G^{*n}(x+\Delta)\ =\
\sum_{n=1}^\infty \theta^n F^{*n}(x+\Delta)
\end{eqnarray*}
is asymptotically eqiuvalent to
\begin{eqnarray*}
&& \sum_{n=1}^\infty n\theta^n F(x+\Delta)\ =\
\frac{\theta}{(1-\theta)^2}F(x+\Delta)
\ =\ \frac{1}{(1-\theta)^2}G(x+\Delta)
\end{eqnarray*}
if and only if $G\in \calS_\Delta$.

Alternatively, one may use the representation
$U=H/(1-\theta)$ where $H$ is the distribution
of $X_1+\cdots+X_\tau$ where
$\P(\tau=n)=(1-\theta)\theta^n$, $n=0,1,2,\ldots$
and the $X_k$ are i.i.d.\ with distribution
$F=G/\theta$ (see [\ref{APQ}] Proposition 2.6
p.\ 114). Hence by Theorem \ref{series}
\begin{eqnarray*}U(x+\Delta)&=&
\frac{H(x+\Delta)}{1-\theta}\ \sim \frac{\E\tau F(x+\Delta)}{1-\theta}
\\ &=&\frac{\theta}{(1-\theta)^2}F(x+\Delta) \ =\
\frac{1}{(1-\theta)^2}G(x+\Delta).
\end{eqnarray*}

We now turn to the renewal equation
\begin{equation}\label{ren}
Z(x)\ =\ z(x)+\int_0^x Z(x-y)\, G(dy),\ \ x\ge 0,
\end{equation}
where $z\ge 0$ and $z$ is locally  bounded. This
together with $\theta\le 1$ is more than sufficient to
ensure that
$$
Z(x)=\int_0^x z(x-y)U(dy)
$$
is the unique locally bounded solution. The key renewal theorem
states that $Z(x)$ has limit $\mu_G^{-1}\int_0^\infty
z(y)dy$  when $\theta=1$. Light-tailed asymptotics
of $Z(x)$ when $\theta<1$ is also available
(see [\ref{F}] or [\ref{APQ}], VI.5) and has found numerous
applications. Therefore, it is surprising that
heavy-tailed asymptotics
when $\theta<1$ appears not to have been discussed before
a specific application came up in Asmussen  [\ref{A2}].
A result was stated there which contains the basic intuition, but
the proof is heuristic as well as the conditions
are not formulated in a precise form.
The analysis of the preceding parts of this paper
allows for a more rigorous treatment, and we shall show
(see [\ref{F}], [\ref{APQ}] for the definition of $z$ to be
d.R.i.\ $=$ directly Riemann integrable):

\begin{Theorem}\label{reneqn}
Assume $\theta<1$ and define
$g(x)=G(x,x+1]$, $I=\int_0^\infty z(y)dy$. Then:

{\rm (i)} if $G\in \calS_\Delta$ for all $T<\infty$,
$z$ is d.R.i., and $z(x)/g(x)\to 0$, then
$$
Z(x)\,\sim\,\frac{I}{(1-\theta)^2}g(x);
$$

{\rm (ii)} if $G\in \calS_\Delta$ for all $T<\infty$,
$z$ is d.R.i., and $z(x)/g(x)\to c\in(0,\infty)$, then
$$
Z(x)\ \sim\
\left(\frac{I}{(1-\theta)^2}+\frac{c}{1-\theta}\right)g(x);
$$

{\rm (iii)} if the probability density $z(y)/I$
is subexponential and $z(x)/g(x)\to \infty$, then
$$
Z(x)\,\sim\,\frac{1}{1-\theta}z(x).
$$
\end{Theorem}

\proof. In (i) and (ii), the assumptions imply $G(x,x+1/n]\sim g(x)/n$
for all $n$ and $g(x+y)/g(x)\to1$ uniformly in $|y|<y_0<\infty$.
Therefore applying Proposition \ref{inter.integral} to
the probability measure $(1-\theta)U$ and
appealing to Proposition \ref{renfct} with $T=1/n$
shows that for each $n$ we can find $h_n(x)\to\infty$ such that
$h_n(x)<x/2$  and
\begin{equation}\label{ren0}
U(x-(k+1)/n,x-k/n]\ \sim\
\frac{g(x)}{n(1-\theta)^2}\ \mbox{ uniformly in }
k\le nh_n(x),
\end{equation}
\begin{eqnarray}\label{ren2}
\int_0^{x-h_n(x)} g(x-y)U(dy) &\sim& (1-\theta)^{-1}g(x),\\
\int_{h_n(x)}^{x-h_n(x)}g(x-y)U(dy) &=& o(g(x))\label{ren1}
\end{eqnarray}
(without loss of generality, we may assume that $nh_n(x)$
is an integer). We will use the decomposition
$Z(x)=$ $J_{1,n}+J_{2,n}+J_{3,n}$ where
$J_{1,n}=$ $\int_0^{h_n(x)}z(x-y)U(dy)$ and similarly
$J_{2,n},J_{3,n}$ are the integrals over
$(h_n(x),x-h_n(x)]$, resp. $(x-h_n(x),x]$.

In (i), we replace $h_n$ by a smaller $h_n$ if necessary
to ensure $z(x-y)/g(x)\to 0$ uniformly in
$|y|\le h_n(x)$ (this is possible since
$g\in\calL$), implying
$J_{1,n}=o(g(x))$. Next,
\begin{eqnarray*}
J_{2,n} &=& o(1)\int_{h_n(x)}^{x-h_n(x)}g(x-y)U(dy)\ =\ o(g(x))
\end{eqnarray*}
by (\ref{ren1}). Finally, writing
$\overline z_n(x)=\sup_{|y-x|\le 1/n}z(y)$, (\ref{ren0}) yields
\begin{eqnarray*}
J_{3,n} &\le&
\sum_{k=0}^{nh_n(x)}\overline z_n(k/n)U(x-(k+1)/n,x-k/n]\\
&\sim&\frac{g(x)/n}{(1-\theta)^2}\sum_{k=0}^{nh_n(x)}\overline z_n(k/n)
\ \sim\ \frac{g(x)/n}{(1-\theta)^2}\sum_{k=0}^\infty\overline z_n(k/n)
\end{eqnarray*}
Since $z$ is d.R.i., $n^{-1}\sum\cdots\to I$ as $n\to\infty$,
yielding $\limsup Z(x)/g(x)\le (1-\theta)^{-2}I$ in (i).
The proof of $\liminf Z(x)/g(x)\ge(1-\theta)^{-2}I$ is similar.

In (ii), we may assume $z(x-y)/g(x)\to c$ uniformly in
$|y|\le h_n(x)$ and then get
$$
J_{1,n}\ \sim\ cg(x)U(h_n(x))\ \sim\ cg(x)U(\infty)\ =\
cg(x)(1-\theta)^{-1}.
$$
For $J_{2,n}$, we have to replace $o(1)$ by $O(1)$,
but the result remains $o(g(x))$. Finally, $J_{3,n}$
can be treated just as in (i), and (ii) is proved.

In (iii), consider the probability measure $K$ with density $z(x)/I$.
The measure $K$ is $\Delta$-subexponential for any $\Delta$.
Put $\Delta=(0,1]$ and write
\begin{eqnarray*}
\int_0^x z(x-y)U(dy) &=&
\int_0^{x-h} z(x-y)U(dy)+\int_{x-h}^x z(x-y)U(dy)
=I_1(x,h)+I_2(x,h),\\
(K*U)(x+\Delta) &=&
\int_0^{x-h} K(x-y+\Delta)U(dy)+\int_{x-h}^{x+1}K(x-y+\Delta)U(dy)\\
&&\hspace{10mm} =I'_1(x,h)+I'_2(x,h).
\end{eqnarray*}
For any fixed $h$ we have, as $x\to\infty$,
\begin{eqnarray*}
I_2(x,h) &\le& h\cdot\sup_{y\le h}|z(y)|\cdot U(x-h,x] = o(z(x))
\end{eqnarray*}
and, by the same reasons, $I'_2(x,h)=o(z(x))$.
Then it is possible to choose $h(x)\uparrow\infty$
such that still $I_2(x,h(x))=o(z(x))$ and $I'_2(x,h(x))=o(z(x))$.
Since $z\in\calL$, $z(x)\sim I\cdot K(x+\Delta)$ and
$I_1(x,h(x))\sim I\cdot I'_1(x,h(x))$.
Combining these estimates we deduce
\begin{eqnarray*}
\int_0^x z(x-y)U(dy) &\sim& I\cdot(K*U)(x+\Delta).
\end{eqnarray*}
Applying Proposition \ref{F.1.F.2} with $G_1=K$, $G_2=U(1-\theta)$,
$c_1=1$, and $c_2=0$ finally yields
\begin{eqnarray*}
I\cdot(K*U)(x+\Delta) &\sim& \frac{I}{1-\theta}K(x+\Delta)
\sim \frac{z(x)}{1-\theta}.
\end{eqnarray*}

\mysection{The compound Poisson process}\label{Poisson}

Let $F$ be a distribution on ${\bf R}_+$ and $\mu$
a positive constant. Let $G$ be
the compound Poisson distribution
\begin{eqnarray*}
G(B) &=& e^{-\mu} \sum_{n\ge0}\frac{\mu^n}{n!} F^{*n}(B).
\end{eqnarray*}

\begin{Theorem}
Let $0<T\le\infty$. If $T<\infty$, then assume
$F\in\calL_\Delta$.
Then the following assertions are equivalent:

{\rm(i)} $F\in\calS_\Delta$;

{\rm(ii)} $G(x+\Delta)\sim \mu F(x+\Delta)$
as $x\to\infty$.
\end{Theorem}

The proof follows from Theorem \ref{series},
with $p_n=\mu^ne^{-\mu}/n!$.

The case $T=\infty$ was considered, for regularly varying
tails, in [\ref{CC}, \ref{C}] and,
for subexponential tails, in [\ref{EGV}].

\mysection{Infinitely divisible laws}\label{IDL}

Let $F$ be an infinitely divisible law on $[0,\infty)$.
The Laplace transform of an infinitely divisible law $F$
can be expressed as
\begin{eqnarray*}
\int_0^\infty e^{-\lambda x}F(dx) &=&
e^{-a\lambda-\int_0^\infty(1-e^{-\lambda x})\nu(dx)}
\end{eqnarray*}
(see, for example, [\ref{F}, p. 450]).
Here $a\ge0$ is a constant and the L\'evy measure $\nu$
is a Borel measure on $(0,\infty)$ with the properties
$\mu=\nu((1,\infty))<\infty$ and $\int_0^1 x\nu(dx)<\infty$.
Put $G(B)=\nu(B\cap(1,\infty))/\mu$.

The relations between the tail behaviour of measure $F$
and the corresponding L\'evy measure $\nu$ were considered
in Theorem 1 in [\ref{EGV}]. We prove the following
local analogue of that result.

\begin{Theorem}
Let $0<T\le\infty$. If $T<\infty$, then assume
$G\in\calL_\Delta$.
Then the following assertions are equivalent:

{\rm(i)} $G\in\calS_\Delta$;

{\rm(ii)} $\nu(x+\Delta)\sim F(x+\Delta)$
as $x\to\infty$.
\end{Theorem}

\proof. It is pointed out in [\ref{EGV}]
that the distribution $F$ admits the representation
$F=F_1*F_2$, where $F_1(x,\infty)=O(e^{-\varepsilon x})$
for some $\varepsilon>0$ and
\begin{eqnarray*}
F_2(B) &=&
e^{-\mu}\sum_{n=0}^\infty \frac{\mu^n}{n!}G^{*n}(B).
\end{eqnarray*}
Now, by Theorem \ref{series}, with
$p_n=\mu^ne^{-\mu}/n!$ we have
$$
F_2(x+\Delta) \sim \mu G(x+\Delta)=\nu(x+\Delta)
\quad\mbox{as }x\to\infty.
$$
Since $F_1(x+\Delta)=o(G(x+\Delta))$ as $x\to\infty$,
by Proposition \ref{F.1.F.2}
$$
F(x+\Delta)=(F_1*F_2)(x+\Delta) \sim F_2(x+\Delta)
\sim\nu(x+\Delta).
$$

\mysection{Branching processes}\label{BP}

In this section we consider the limit behaviour of
sub-critical, age-dependent branching processes for which
the Malthusian parameter does not exist.

Let $h(z)$ be the particle production generating function
of an age-dependent branching process with particle lifetime
distribution $F$ (see [\ref{AN}, Chapter IV], [\ref{H}, Chapter VI]
for background). We take the process to be sub-critical,
i.e. $A\equiv h'(1)<1$. Let $Z(t)$ denote the number
of particles at time $t$. It is known
(see, for example, [\ref{AN}, Chapter IV, Section 5] or
[\ref{Ch}]) that $A(t)={\bf E}Z(t)$ admits the representation
\begin{eqnarray}\label{A.represent}
A(t) &=& (1-A)\sum_{n=1}^\infty A^{n-1}(1-F^{*n}(t)).
\end{eqnarray}
It was proved in [\ref{Ch}] for sufficiently small values of $A$
and then in [\ref{CNW}, \ref{CNW2}] for any $A<1$ that
$A(t)\sim\overline F(t)/(1-A)$ as $t\to\infty$, provided $F\in\calS$.

Applying Theorem \ref{series} with $p_n=(1-A)A^{n-1}$
(see also Proposition \ref{renfct}), we deduce

\begin{Theorem}
If $F\in\calL_\Delta$, then the following are equivalent:

{\rm(i)} $F\in\calS_\Delta$;

{\rm(ii)} $A(t)-A(t+T)\sim F(t+\Delta)/(1-A)$ as $t\to\infty$.
\end{Theorem}

\begin{center}
\bf Acknowledgment
\end{center}

The authors gratefully acknowledge
helpful discussions with Stan Zachary.

\section*{\centering\small References}

\newcounter{bibcoun}
\begin{list}{\arabic{bibcoun}.}{\usecounter{bibcoun}\itemsep=0pt}
\small

\item\label{APQ}
   Asmussen, S. (1987).
   {\it Applied Probability and Queues},
   Wiley, Chichester (2nd ed.\ (2001) to  be published by
   Springer, New York).

\item\label{A}
   Asmussen, S. (2000).
   {\it Ruin Probabilities},
   World Scientific, Singapore.

\item\label{A2}
   Asmussen, S. (1998).
   A probabilistic look at the Wiener--Hopf equation.
   {\it SIAM Review} {\bf 40}, 189--201.

\item\label{AKKKT}
   Asmussen, S., Kalashnikov, V., Konstantinides, D.,
   Kl\"uppelberg, C., and Tsitsiashvili, G. (2002).
   A local limit theorem for random walk maxima
   with heavy tails.
   {\it Statist.\ Probab.\ Letters} {\bf 56}, 399--404.

\item\label{AN}
   Athreya, K., and Ney, P. (1972).
   {\it Branching Processes},
   Springer, Berlin.

\item\label{BD}
   Bertoin, J., Doney, R. A. (1994).
   On the local behaviour of ladder height distributions.
   {\it J.\ Appl.\ Prob.} {\bf 31}, 816--821.

\item\label{B}
   Borovkov, A. A. (1976).
   {\it Stochastic Processes in Queueing Theory},
   Springer, Berlin.

\item\label{CC}
   Callaert, H., and Cohen, J. W. (1972).
   A lemma on regular variation of a transient renewal function.
   {\it Z. Wahrscheinlichkeitstheorie verw. Gebiete}
   {\bf 24}, 275--278.

\item\label{Ch}
   Chistyakov, V. P. (1964).
   A theorem on sums of independent positive random variables
   and its application to branching random processes,
   {\it Theor. Probability Appl.} {\bf 9}, 640--648.

\item\label{CNW}
   Chover, J., Ney, P., and Wainger, S. (1973).
   Functions of probability measures.
   {\it J. d'Analyse Math\'ematique} {\bf 26}, 255--302.

\item\label{CNW2}
   Chover J., Ney P., and Wainger S. (1973).
   Degeneracy properties of subcritical branching processes.
   {\it Ann. Probab.} {\bf 1}, 663--673.

\item\label{C}
   Cohen, J. W. (1973).
   Some results on regular variation for distributions
   in queueing and fluctuation theory.
   {\it J. Appl. Probab.} {\bf 10}, 343--353.

\item\label{EKM}
   Embrechts, P., Kl\"uppelberg, C., and Mikosch, T. (1997).
   {\it Modelling Extremal Events for Insurance and Finance},
   Springer, Berlin.

\item\label{EG}
   Embrechts, P. and Goldie, C. M. (1982).
   On convolution tails.
   {\it Stoch. Processes Appl.} {\bf 13}, 263--278.

\item\label{EGV}
   Embrechts, P., Goldie, C. M., and Veraverbeke, N. (1979).
   Subexponentiality and infinite divisibility.
   {\it Z. Wahrscheinlichkeitstheorie verw. Gebiete}
   {\bf 49}, 335--347.

\item\label{EV}
   Embrechts, P., and Veraverbeke, N. (1982).
   Estimates for the probability of ruin with special emphasis
   on the possibility of large claims.
   {\it Insurance: Math. and Economics} {\bf 1}, 55--72.

\item\label{F}
   Feller, W. (1971).
   {\it An Introduction to Probability Theory and Its
   Applications} Vol. 2,
   Wiley, New York.

\item\label{FZ}
   Foss, S., and Zachary, S. (2002).
   The maximum on a random time interval of a random walk
   with long-tailed increments and negative drift.
   To appear in {\it Ann.\ Appl.\ Prob.}

\item\label{H}
   Harris, T. (1963).
   {\it The Theory of Branching Processes},
   Springer, Berlin.

\item\label{Claudia}
   Kl\"uppelberg, C. (1989).
   Subexponential distributions and characterization of
   related classes.
   {\it Probab. Th. Rel. Fields} {\bf 82}, 259--269.

\item\label{K}
   Korshunov, D. (1997).
   On distribution tail of the maximum of a random walk.
   {\it Stoch. Proc. Appl.} {\bf 72}, 97--103.

\item\label{RSST}
   Rolski, T., Schmidli, H., Schmidt, V, and Teugels, J.
   (1998).
   {\it Stochastic Processes for Insurance and Finance},
   Wiley, Chichester.

\item\label{S}
    Sgibnev, M. S.
    Banach algebras of functions that have identical
    asymptotic behaviour at infinity.
    {\it Siberian Math. J.} {\bf 22}, 179--187.

\item\label{T}
    Teugels, J. L. (1975).
    The class of subexponential distributions.
    {\it Ann. Probab.} {\bf 3}, 1000--1011.

\item\label{Ver}
    Veraverbeke, N. (1977).
    Asymptotic behavior of Wiener-Hopf factors
    of a random walk.
    {\it Stoch. Proc. Appl.} {\bf 5}, 27--37.

\end{list}

\end{document}